\documentclass[11pt]{amsart}

\usepackage[colorlinks=true, pdfstartview=FitV, linkcolor=blue, 
citecolor=blue]{hyperref}

\usepackage{amssymb,amsmath,accents}
\usepackage{bbm}
\usepackage{graphicx}
\usepackage{a4wide}

\theoremstyle{definition}

\newcommand{\ts}{\hspace{0.5pt}}
\newcommand{\nts}{\hspace{-0.5pt}}

\newcommand{\RR}{\mathbb{R}\ts}
\newcommand{\CC}{\mathbb{C}}
\newcommand{\ZZ}{\mathbb{Z}}

\newcommand{\NN}{\mathbb{N}}

\newcommand{\QQ}{\mathbb{Q}}

\newcommand{\cB}{\mathcal{B}}

\newcommand{\cL}{\mathcal{L}}

\newcommand{\vL}{\varLambda}

\newcommand{\ee}{\ts\mathrm{e}}
\newcommand{\ii}{\mathrm{i}\ts}

\DeclareMathOperator{\dens}{dens}

\DeclareMathOperator{\vol}{vol}
\DeclareMathOperator{\real}{Re}
\DeclareMathOperator{\imag}{Im}
\DeclareMathOperator{\im}{im}

\newcommand{\widec}{\smash{\raisebox{16.5pt}{\rotatebox{180}
{$\widehat{\hskip 12pt}$}}}}
\newcommand{\widecl}{\smash{\raisebox{17.5pt}{\rotatebox{180}
{$\widehat{\hskip 24pt}$}}}}
\newcommand{\widecheck}[1]{\hskip 2pt\widec \hskip -14pt {#1}}
\newcommand{\widecheckl}[1]{\hskip 2pt\widecl \hskip -26pt {#1}}

\newcommand{\sB}{\ts\underline{\nts B\!}\,}
\newcommand{\sR}{\ts\underline{\nts R\nts\nts}\ts\ts}
\newcommand{\sQ}{\ts\underline{\nts Q \nts}\ts}
\newcommand{\defeq}{\mathrel{\mathop:}=}

\begin{document}

\title{Diffraction of a model set with complex windows}

\author{Michael Baake}
\address{Fakult\"at f\"ur Mathematik, Universit\"at Bielefeld, \newline
       \indent  Postfach 100131, 33501 Bielefeld, Germany}
\email{mbaake@math.uni-bielefeld.de}

\author{Uwe Grimm}
\address{School of Mathematics and Statistics,
  The Open University,\newline 
  \indent Walton Hall, Milton Keynes MK7 6AA, UK} 
\email{uwe.grimm@open.ac.uk}

\begin{abstract} 
The well-known plastic number substitution gives 
rise to a ternary inflation tiling of
the real line whose inflation factor is the smallest
Pisot--Vijayaraghavan number. 
The corresponding dynamical system has pure point
spectrum, and the associated control point sets
can be described as regular model sets whose
windows in two-dimensional internal space 
are Rauzy fractals with a complicated structure.
Here, we calculate the resulting pure point 
diffraction measure via a Fourier matrix
cocycle, which admits a closed formula for 
the Fourier transform of the Rauzy fractals,
via a rapidly converging infinite product.
\end{abstract}

\maketitle

Consider the primitive substitution $\varrho$ defined by
$a\mapsto b\mapsto c \mapsto ab$ on the ternary
alphabet $\{ a,b,c\}$; see \cite[Ex.~4.4]{TAO} for details.  
Its substitution matrix reads
\[
  M \, = \, \begin{pmatrix} 0 & 0 & 1\\
    1 & 0 & 1\\ 0 & 1 & 0\end{pmatrix} ,
\]
with characteristic polynomial $x^3-x-1$. The latter is
irreducible, with one real root, 
\[
  \beta\, = \, \frac{(9+\sqrt{69})^{\frac{1}{3}} +
  (9-\sqrt{69})^{\frac{1}{3}}}{18^{\frac{1}{3}}}
   \, \approx \, 1.32472 \ts ,
\]
and a complex conjugate pair $\alpha, \overline{\alpha}$ of algebraic
conjugates, where $\lvert\alpha\rvert^2=\beta^{-1}=\beta^2-1$. Below,
we use $\alpha$ for the number with positive imaginary part.

The algebraic integer $\beta$ is a unit, and the smallest
Pisot--Vijayaraghavan (PV) number, known as the \emph{plastic number};
compare \cite[p.~50 and Ex.~2.17]{TAO}.  It is also the
Perron--Frobenius eigenvalue of $M$, where the corresponding left and
right eigenvectors read
\[
  \langle u\ts | \, = \,
  \frac{3+\beta+7\beta^2}{23}(1,\beta,\beta^2)
  \quad \text{and} \quad
    |\ts v\rangle \, = \,
    (2-\beta^2,\beta^2-\beta,\beta-1)^T ,
\]
which are strictly positive, and normalised such that
$\langle 1 | \ts v \rangle = \langle u \ts | \ts v \rangle =1$. So,
the entries of $|\ts v\rangle$ are the relative frequencies of the
letters in the symbolic sequences defined by $\varrho$, and
$P \defeq |\ts v\rangle \langle u \ts |$ is a projector of rank $1$,
that is, $P^2 = P$ and $\im (P) = \RR \ts |\ts v\rangle$.

To turn the symbolic sequences into tilings, we choose intervals of
natural length, namely $1,\beta$ and $\beta^2$, 
with control points on their left endpoints. 
For each sequence in the symbolic hull of $\varrho$, this
leads to a typed point set, $\vL = \vL_{a} \ts\ts\dot{\cup}\ts\ts \vL_{b} 
\ts\ts\dot{\cup}\ts\ts \vL_{c}$, from the three types of control points. The 
average distance between neighbouring points in $\vL$ is
$\bar{s}=4 + 2 \beta - 3 \beta^2$, so we get
$\dens(\vL) = 1/\bar{s} = \frac{1}{23} (3+\beta+7 \beta^2)$ 
and $\dens (\vL^{}_{i}) = v^{}_{i} \ts \dens (\vL)$ for $i \in \{ a,b,c \}$.

To continue, we work with the rank-$3$ $\ZZ$-module
$L=\ZZ[\beta]=\langle 1,\beta,\beta^2\rangle^{}_{\ZZ}$, 
which comprises all possible coordinates of our control points
(relative to one of them, then considered as $0$).
The Minkowski embedding \cite[Ex.~3.5]{TAO} of $L$
is a lattice $\cL$ in $\RR \!\times\! \CC \simeq \RR^3$,
with $\dens (\cL) = 2/\sqrt{23}$. A canonical choice for the 
basis matrix of $\cL$ and its dual, $\cL^*$, is
\[
    \cB\, = \, \begin{pmatrix}
  1 & \beta & \beta^2 \\
  1 & \real(\alpha) & \real(\alpha^2) \\
  0 & \imag(\alpha) & \imag(\alpha^2)
  \end{pmatrix} \quad\text{and}\quad
    \cB^{*} \, = \, \frac{2}{\sqrt{23}} \begin{pmatrix}
  \imag(\alpha) \beta^{-1} &
  \imag(\alpha) \beta &
  \imag(\alpha) \\
  2\imag(\alpha)\beta^2 & 
  -\imag(\alpha)\beta &
  -\imag(\alpha)\\
  \beta & -1+\frac{3}{2}\beta^2 & -\frac{3}{2}\beta
\end{pmatrix},
\]
where $\real(\alpha) = -\beta/2$, $\real(\alpha^2)=1-\beta^2/2$,
$\imag(\alpha) = \frac{4+9\beta-6\beta^2}{2\sqrt{23}}$ and
$\imag(\alpha^2) = \frac{6+2\beta-9\beta^2}{2\sqrt{23}}$
in our setting. We note in passing that duality can
also be defined via the quadratic form
$(x,y) \defeq xy + \real\bigl(\sigma(x)\ts
\overline{\sigma(y)}\ts\bigr)$, where
$\sigma \! : \, \QQ (\beta) \longrightarrow \QQ (\alpha)$ denotes the
algebraic conjugation map defined by $\beta\mapsto\alpha$
and its unique extension to a field isomorphism.

We can extract the  \emph{Fourier module} from the 
first line of the dual basis matrix $\cB^*$ as
\begin{equation}\label{eq:F-mod}
   L^{\circledast} \,=\, \frac{2\imag(\alpha)}{\beta\sqrt{23}}
   \, \langle 1,\beta^2,\beta\rangle^{}_{\ZZ} \,=\,
   \frac{5-6\beta+4\beta^2}{23} \ts L \ts ,
\end{equation}
which is also the dynamical spectrum (in additive notation) of the
tiling dynamical system as well as the (topologically equivalent)
model set dynamical system defined by the control point sets; see
\cite{BL-review} for background.  Our system has pure point
diffraction spectrum, because the defining point set
is a regular model set (see below).

The Bragg peaks can be indexed by three integer Miller indices
$(n^{}_{0},n^{}_{1},n^{}_{2})$, where we parameterise the wave number,
in line with \eqref{eq:F-mod}, as
\begin{equation}\label{eq:wave-number}
  k \, = \, k (n^{}_{1}, n^{}_{2}, n^{}_{3} ) \, = \,
  \frac{5-6\beta+4\beta^2}{23} \bigl(n^{}_{0}+
  n^{}_{1}\beta +n^{}_{2}\beta^2 \bigr)  .
\end{equation}  
The original version of the Minkowski embedding of $\ZZ [ \beta]$
uses complex numbers, which is natural from an algebraic
perspective. However, for Fourier analysis, we better work with 
real numbers, via $\CC\simeq \RR^2$, as initiated above
for $\cB$ and $\cB^*$. Given the conjugation map $\sigma$, the
$\star$-map: $\QQ(\beta)\longrightarrow \RR^2$ is conveniently defined
by $x \mapsto \bigl(\real(\sigma(x)),\imag(\sigma(x))\bigr)^T$.
If $k\in L^{\circledast}$ is parameterised as in \eqref{eq:wave-number},
this means
\[
  k^{\star}\, =\, \begin{pmatrix}
   \frac{1}{46}\bigl((18 n^{}_{0} - 4 n^{}_{1} + 6 n^{}_{2})
    + (6 n^{}_{0} - 9 n^{}_{1} + 2 n^{}_{2}) \beta 
    - (4 n^{}_{0} - 6 n^{}_{1} + 9 n^{}_{2}) \beta^2 \bigr) \\
    \frac{1}{2\sqrt{23}} \bigl( 2 n^{}_{2} + 
    (3 n^{}_{1} - 2 n^{}_{0}) \beta - 3 n^{}_{2} \beta^2\bigr)
   \end{pmatrix} .
\]

Next, we need the displacement (or digit) matrix $T$ 
of our inflation rule in direct space; see~\cite{BFGR} for 
definitions and properties.
With the interval lengths as chosen above, it reads
\[
    T\, = \, \begin{pmatrix} \varnothing & \varnothing & \{0\}\\
    \{0\} & \varnothing & \{1\}\\
    \varnothing & \{0\} & \varnothing
    \end{pmatrix} ,
\]
which is also reflected in the point set iteration
\[
   \vL_{a}^{\prime} \, = \, \beta\vL^{}_{c}\, , \quad
   \vL_{b}^{\prime} \, = \, \beta\vL^{}_{a} \cup 
    (\beta\vL^{}_{c}+1)\, ,\quad
   \vL_{c}^{\prime} \, = \, \beta\vL^{}_{b} \ts .
\]
It is known that we obtain a regular model set \cite{Bernd}, 
with three topologically regular windows that have disjoint interiors. 
They are compact sets $W_{\! a}, W_{\nts b}, W_{\! c} 
\subset \CC$ that satisfy
\begin{equation}\label{eq:win-IFS}
    W_{\! a}^{} \, = \, \alpha W^{}_{\! c}\, ,\quad
    W_{\nts b}^{} \, = \, \alpha W^{}_{\! a} \cup 
    (\alpha W^{}_{\! c}+1)\, ,\quad
    W_{\! c}^{} \, = \, \alpha W^{}_{\! b} \ts ,
\end{equation}
which defines a contractive iterated function system 
(IFS), so that the fixed point (and hence our
window triple) is unique. Due to the nature of $\vL$
and $\vL_i$ as regular model sets, the volumes (areas) 
of the windows, compare \cite[Ex.~7.6 and Fig.~7.3]{TAO},
satisfy the relations $\dens (\cL) \vol (W^{}_i) = v^{}_{i} \ts
\dens (\vL)$, which results in the values
\[
     \imag (\alpha) \bigl( \beta^2 - 1, \beta , 1 \bigr)
     \, = \, \tfrac{1_{\phantom{j}}}{2 \sqrt{23}} \,
     \bigl( 5-6 \beta + 4 \beta^2 , 
     -6 -2 \beta + 9 \beta^2 , 4 + 9 \beta - 6 \beta^2
     \bigr).
\]

\begin{figure}
\centerline{\includegraphics[width=0.6\textwidth]{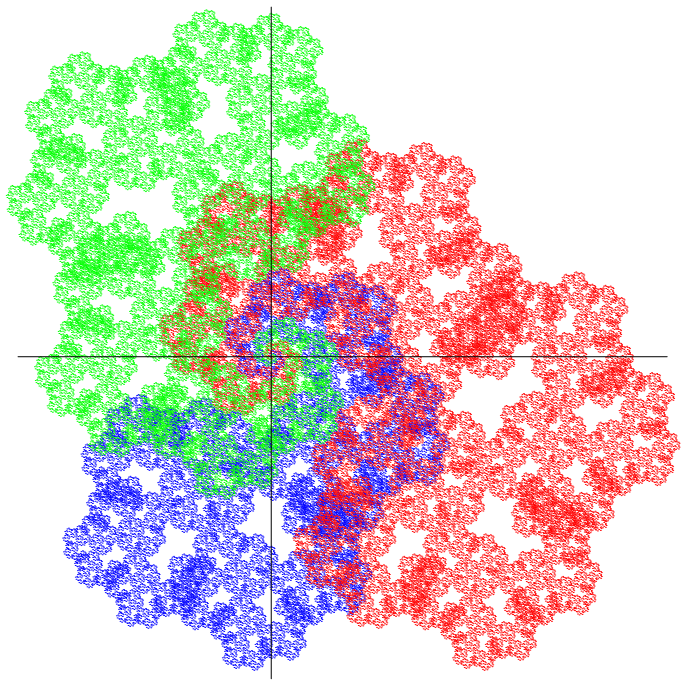}}
\caption{\label{winfig}Illustration of the three windows $W_{\! a}$
  (blue), $W_{\nts b}$ (red) and $W_{\! c}$ (green) of the plastic
  number inflation rule. For a more detailed version of 
  $W_{\nts b}$ and further properties, see \cite[Fig.~7.3]{TAO}.}
\end{figure}

The three windows are complex Rauzy fractals with a complicated
topological structure \cite{ST}; see Figure~\ref{winfig} for an
illustration and \cite[Sec.~7.4]{PFBook} and references therein for
background on Rauzy fractals. From the general diffraction result for
regular model sets, see \cite{TAO} and references therein, the
weighted Dirac comb $\omega = \sum_{i\in\{a,b,c\}}h_{i}\,
\delta^{}_{\!\vL_{i}}$ with $h_{i}\in\CC$ has diffraction
\begin{equation}\label{eq:intens}
    \widehat{\gamma^{}_{\omega}}\, = 
    \sum_{k\in L^{\circledast}} I(k)\, \delta^{}_{k}\quad\text{with}\quad
   I(k) \, = \, \Bigl| \sum_{i} h_{i} \, A_{i}(k) \Bigr|^{2} ,
\end{equation}
because the Bombieri{\ts}--Taylor property holds for primitive inflation
rules \cite{BGM}.  Here, setting $W=W_{\! a}\cup W_{\nts b}\cup W_{\!
  c}$, the amplitudes are given by
\[
     A_{i}(k) \, = \, \frac{\dens(\vL_{i})}{\vol(W_i)} \, 
     \widecheck{1^{}_{W_i}}(k^{\star}) \, = \, 
     \frac{\dens(\vL)}{\vol(W)} \, 
     \widecheck{1^{}_{W_i}}(k^{\star}) \ts ,
\]
where $1^{}_{K}$ denotes the characteristic function of the set $K$.
Though this formula looks nice, it is difficult to calculate
$\widecheck{1^{}_{W_i}}$ directly, due to the fractal nature of the
windows.  Let us next explain an alternative approach that harvests
the inflation nature of our point sets.

We start by reformulating the window IFS \eqref{eq:win-IFS} in
$\RR^2$. Recall that multiplication with $\alpha\in\CC$ is matrix
multiplication with $\sQ=\left(\begin{smallmatrix} \real(\alpha) &
  -\imag(\alpha)\\ \imag(\alpha) &
  \real(\alpha) \end{smallmatrix}\right)$ in $\RR^2$, where
$\det (\sQ) = \lvert \alpha \rvert^2 = \beta^{-1}$.  As the windows
are interior-disjoint, Eq.~\eqref{eq:win-IFS}, 
turns into equations in $L^{1} (\RR^2)$
for the characteristic functions of the windows
and their images. Setting $f^{}_{i} (y) \defeq
\widecheck{1^{}_{W_{i}}} (y)$ and using the relation
\[
    \widecheckl{1^{}_{\! A K \nts + t}} (y) \, = \,
    \lvert\det (A)\rvert \, \ee^{2 \pi \ii \langle t | y \rangle}
    \,\widecheck{1^{}_{\nts K}} (A^T y) \ts ,
\]
which holds for any real, invertible $2\!\times\! 2$-matrix 
$A$ and compact set $K\subset \RR^2$, one finds
\begin{equation}\label{recursion}
     \begin{pmatrix} f^{}_a \\ f^{}_b \\ f^{}_c \end{pmatrix}
     (y) \, = \, \beta^{-1} \sB (y)
     \begin{pmatrix} f^{}_a \\ f^{}_b \\ f^{}_c \end{pmatrix}
     (\sR \ts y)  \ts , \quad \text{with } \, \sR = \sQ^T
     \; \text{and } \,
       \sB(y)\, = \, \begin{pmatrix}
    0 & 0 & 1 \\
    1 & 0 & \ee^{2\pi\ii y^{}_{1}}\\
    0 & 1 & 0 \end{pmatrix} .
\end{equation}
Here, $\sB (y)$ is the internal Fourier matrix that emerges
from the (inverse) Fourier transform of the $\star$-image of
$T$, that is, $\sB_{ij} (y) = \sum_{x\in T_{ij}} 
\ee^{2 \pi \ii \langle x^{\star} | \ts y \rangle}$, where
$y\in\RR^2$.

\begin{figure}
\centerline{\includegraphics[width=0.8\textwidth]{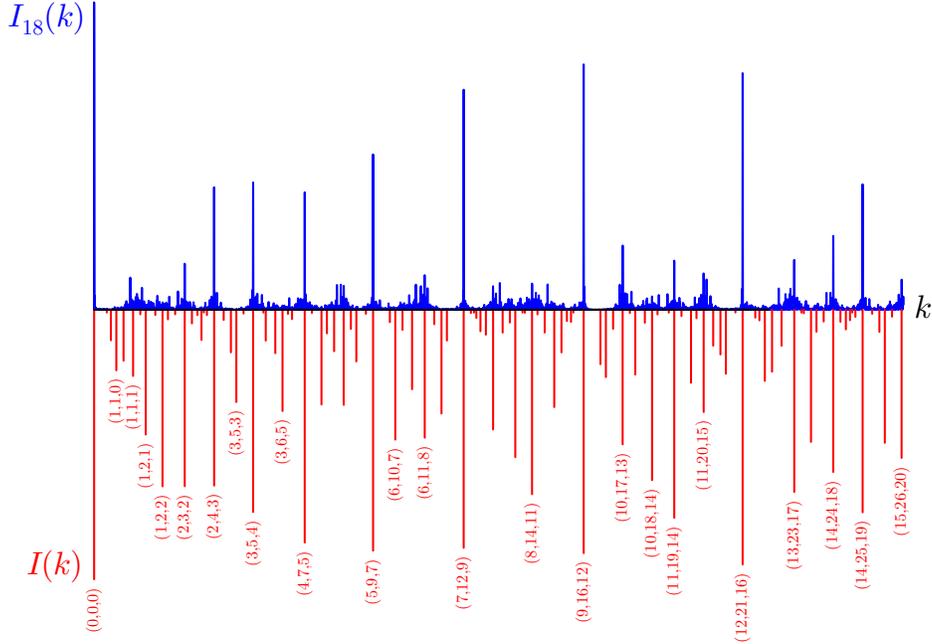}}
\caption{\label{difffig}Bragg peaks of the plastic number inflation
  rule (red, bottom; some with their Miller index triple) as 
  obtained from the cocycle approach, in
  comparison with a finite-size approximation by exponential sums
  (blue, top).}
\end{figure}

Now, for $n\in\NN$, we can define a cocycle via
$\sB^{(n)}(y) = \sB(y) \sB(\sR \ts y) \cdots \sB(\sR^{n-1} y)$, with 
$\sB^{(1)} = \sB$ and $\sB^{(n)}(0)=M^n$. Next,
consider the matrix function
$C(y)\defeq \lim_{n\to\infty} \beta^{-n} \sB^{(n)}(y)$,
which is well defined because the limit exists for
every $y\in\RR^2$. In fact, $C$ is a continuous
matrix function, and one has
\[
   C(0) \, = \, P \, = \, |\ts v\rangle\langle u\ts | \, = \, 
\frac{3+\beta+7\beta^2}{23} \begin{pmatrix}
   2-\beta^2 & \beta - 1 & \beta^2 - \beta \\
   \beta^2 - \beta & 1 + \beta - \beta^2 & \beta^2 - 1 \\
   \beta - 1 & \beta^2 - \beta & 1 + \beta - \beta^2
\end{pmatrix} .
\]
Moreover, one can show that
$C(y) = |\ts c(y)\rangle\langle u\ts |$ with 
$|\ts c(0)\rangle = |\ts v\rangle$ and
$| f (y) \rangle = \frac{\dens (\vL)}{\dens (\cL)} \, | \ts c(y) \rangle$.
In particular, one has $|\ts c(y) \rangle = C(y) \ts
|\ts v \rangle$, which makes the functions $f_i$
accessible.

For $k\in L^{\circledast}$ and $i \in \{ a,b,c \}$, our amplitudes are
\[
    A^{}_{i} (k) \, = \, \dens (\vL) \, c^{}_{i} (k^{\star}) 
    \, = \, \frac{3+\beta+7\beta^2}{23} \, c^{}_{i}(k^{\star}) \ts ,
\]
and the corresponding intensities follow from Eq.~\eqref{eq:intens}.
For any index triple, the intensity at the corresponding wave number
$k\in L^{\circledast}$ is approximated by truncating the infinite
product representation for $C(k^{\star})$ and calculating the
amplitudes as explained above. Here, by using about 50--100 terms, one
obtains the peaks with a relative precision better than $10^{-10}$.

\begin{figure}
\centerline{\includegraphics[width=0.8\textwidth]{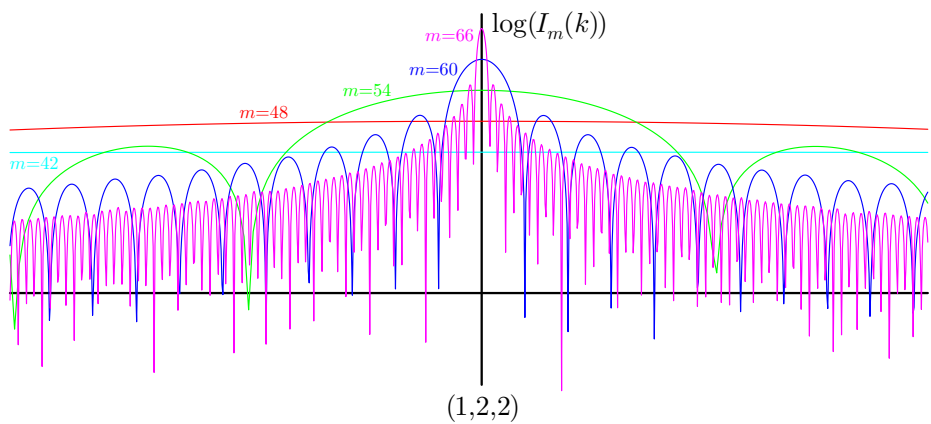}}
\caption{\label{compare}Example of a Bragg peak and its successive
  approximation, for wave numbers $k\in [1.2672395,1.2672405]$ and finite systems
  defined by the point sets from the inflation words $\varrho^{m}(a)$
  with $m\in\{42, 48, 54, 60, 66\}$. The vertical black line denotes
  the location of the Bragg peak $(1,2,2)$. Note that, to be able to
  compare different system sizes, the logarithm of the intensity is
  plotted.}
\end{figure}

Since the diffraction measure $\widehat{\gamma^{}_{\omega}}$ can also
be approximated (in the vague topology) by the absolute squares of
exponential sums of finite approximants (divided by the system size),
we illustrate the diffraction of the uniform Dirac comb (all $h_i=1$)
in Figure~\ref{difffig} in comparison to an approximation for a finite
system with $114$ tiles, obtained by $m=18$ inflation steps from an
initial tile of type $a$. Here, as well as in Figure~\ref{compare},
the function $I_{m}(k)$ for the approximation refers to the (absolutely
continuous) diffraction intensity for the corresponding finite
system. The apparent mismatch emerges from the extremely slow
convergence of the finite-size approximation, which should not come as
a surprise due to the complicated fractal nature of the windows.  To
expand on the latter, we compare the neighbourhood of the peak with
Miller indices $(1,2,2)$ at $k=(1+8\beta+10\beta^2)/23 \approx
1.267240014$ with approximations for various system sizes in
Figure~\ref{compare}.

\begin{figure}
\centerline{\includegraphics[width=0.4\textwidth]{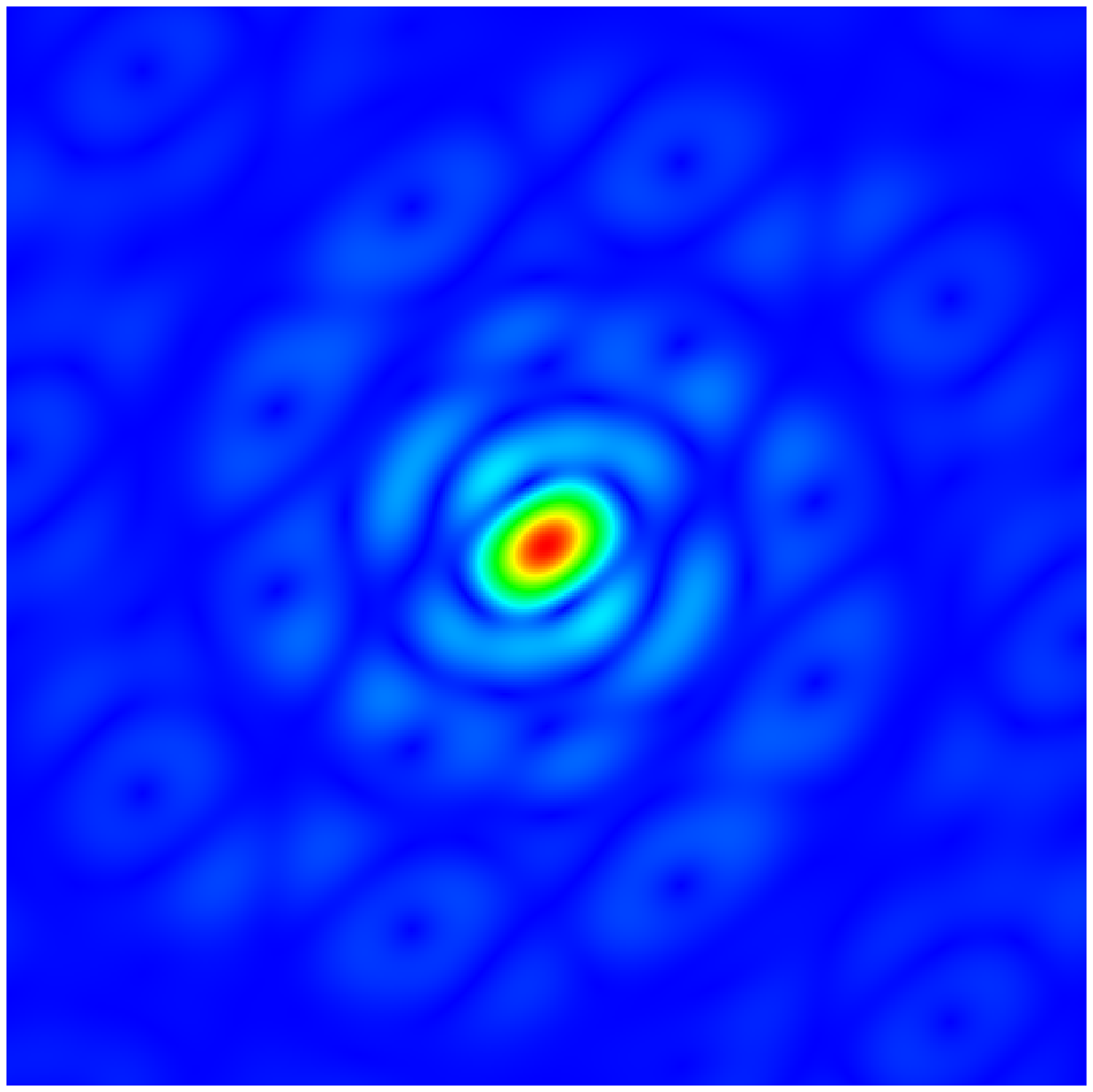}
\hspace{0.05\textwidth}\includegraphics[width=0.4\textwidth]{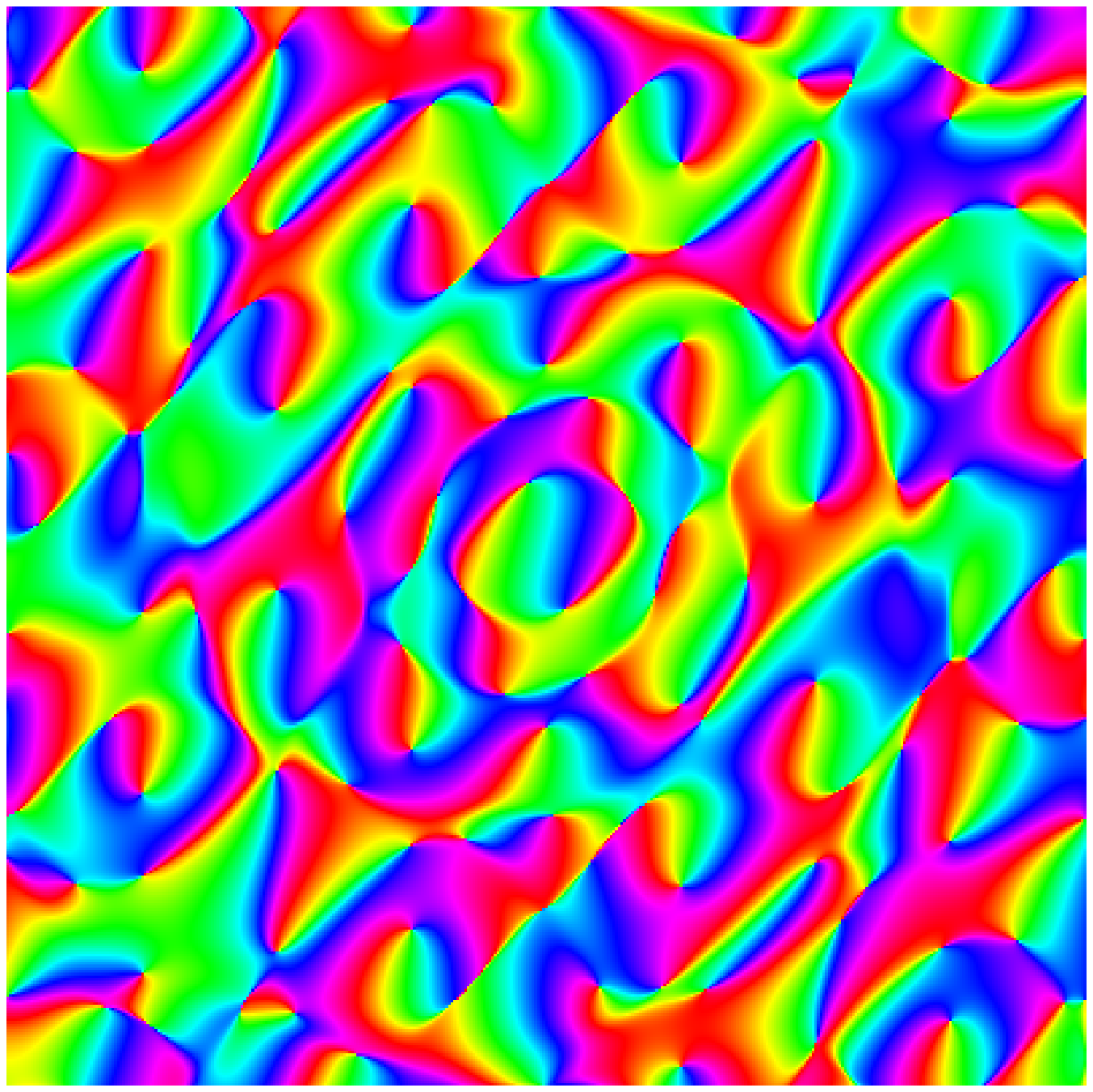}}
\caption{\label{fouriertrans}Inverse Fourier transform $f_{b}(y)$ of
  the window $W_{\nts b}$ for $y \in [-4,4]^2$. The left panel shows
  the absolute value, which takes values between $0$ (blue) and
  $f_{b}(0)=\vol(W_{\nts b})= v^{}_{2}=\beta^2-\beta$ (red). The right
  panel shows the corresponding argument (with red corresponding to
  phase $0$).}
\end{figure}

While the computation of the spectrum requires that $k\in
L^{\circledast}$, the recursion equations \eqref{recursion} for the
(inverse) Fourier transforms of the three windows can be used for any
$y\in\mathbb{R}^{2}$. Thus, the (inverse) Fourier transform of the
windows can be computed efficiently, despite the complex nature of the
windows.  As an example, the inverse Fourier transform
$f_{b}=\widecheck{1^{}_{W_{\nts b}}}$ of the largest window is shown
in Figure~\ref{fouriertrans}. In general, it is difficult to compute
the Fourier transform of fractals; see, for instance, \cite{DF} 
for a related result on self-similar fractals of zero Lebesgue measure.

The cocycle method explained for this example works in full generality
for all primitive PV inflation rules (with a unit inflation factor)
that lead to regular model sets.  In fact, it can even be applied to
any primitive PV inflation rule, as well as to $S$-adic inflations
with the same PV inflation multiplier and compatible tile sizes. When
the spectrum is mixed, the cocycle method reproduces the pure point
part; further details will be given in a forthcoming publication.

\section*{Acknowledgements}

It is our pleasure to thank the MFO at Oberwolfach  
for hospitality, where this manuscript was completed.  
Our work was supported by the German Research
Foundation (DFG), within the CRC 1283 at Bielefeld University, 
and by EPSRC through grant EP/S010335/1.

\end{document}